
\input amssym.def
\input amssym

\documentstyle{amsppt}
\magnification=1200

\font\cal=cmsy10

\hoffset 0,7truecm
\voffset 1truecm
\hsize 15,2truecm
\vsize 22truecm 

\font\eufrak=eufm10

\catcode`\@=11
\long\def\thanks#1\endthanks{%
  \ifx\thethanks@\empty \gdef\thethanks@{%
     \tenpoint#1} 
     \else \expandafter\gdef\expandafter\thethanks@\expandafter{%
     \thethanks@\endgraf#1}%
  \fi}
\def\keywords{\let\savedef@\keywords
  \def\keywords##1\endkeywords{\let\keywords\savedef@
  \toks@{\def\usualspace{{\it\enspace}}\tenpoint} 
  \toks@@{##1\unskip.}%
  \edef\thekeywords@{\the\toks@\frills@{%
  {\noexpand\it 
  Key words and phrases.\noexpand\enspace}}\the\toks@@}}%
  \nofrillscheck\keywords}
\def\makefootnote@#1#2{\insert\footins
 {\interlinepenalty\interfootnotelinepenalty
 \tenpoint\splittopskip\ht\strutbox\splitmaxdepth\dp\strutbox 
 \floatingpenalty\@MM\leftskip\z@\rightskip\z@
 \spaceskip\z@\xspaceskip\z@
 \leavevmode{#1}\footstrut\ignorespaces#2\unskip\lower\dp\strutbox
 \vbox to\dp\strutbox{}}}
\def\subjclass{\let\savedef@\subjclass
  \def\subjclass##1\endsubjclass{\let\subjclass\savedef@
  \toks@{\def\usualspace{{\rm\enspace}}\tenpoint}
  \toks@@{##1\unskip.}%
  \edef\thesubjclass@{\the\toks@\frills@{{%
  \noexpand\rm1991 {\noexpand\it Mathematics
  Subject Classification}.\noexpand\enspace}}%
  \the\toks@@}}%
\nofrillscheck\subjclass}
\def\abstract{\let\savedef@\abstract
 \def\abstract{\let\abstract\savedef@
  \setbox\abstractbox@\vbox\bgroup\noindent$$\vbox\bgroup
  \def\envir@end{\endabstract}\advance\hsize-2\indenti
  \def\usualspace{\enspace}\tenpoint \noindent 
  \frills@{{\smc Abstract.\enspace}}}%
 \nofrillscheck\abstract}

\outer\def\endtopmatter{\add@missing\endabstract
 \edef\next{\the\leftheadtoks}\ifx\next\empty
  \expandafter\leftheadtext\expandafter{\the\rightheadtoks}\fi
 \ifmonograph@\else
   \ifx\thesubjclass@\empty\else \makefootnote@{}%
        {\thesubjclass@}\fi
   \ifx\thekeywords@\empty\else \makefootnote@{}%
        {\thekeywords@}\fi  
   \ifx\thethanks@\empty\else \makefootnote@{}%
        {\thethanks@}\fi   
 \fi
  \pretitle
  \begingroup 
  \ifmonograph@ \topskip7pc \else \topskip4pc \fi
  \box\titlebox@
  \endgroup
  \preauthor
  \ifvoid\authorbox@\else \vskip2.5pcplus1pc\unvbox\authorbox@\fi
  \preaffil
  \ifvoid\affilbox@\else \vskip1pcplus.5pc\unvbox\affilbox@\fi
  \predate
  \ifx\thedate@\empty\else
       \vskip1pcplus.5pc\line{\hfil\thedate@\hfil}\fi
  \preabstract
  \ifvoid\abstractbox@\else
       \vskip1.5pcplus.5pc\unvbox\abstractbox@ \fi
  \ifvoid\tocbox@\else\vskip1.5pcplus.5pc\unvbox\tocbox@\fi
  \prepaper
  \vskip2pcplus1pc\relax}
\def\foliofont@{\tenrm} 
\def\headlinefont@{\tenpoint} 


\def\refstyle#1{\uppercase{%
  \if#1A\relax \def\keyformat##1{[##1]\enspace\hfil}%
  \else\if#1B\relax 
      \def\keyformat##1{\aftergroup\kern
           aftergroup-\aftergroup\refindentwd}%
      \refindentwd\parindent 
  \else\if#1C\relax 
      \def\keyformat##1{\hfil##1.\enspace}%
  \else\if#1D\relax 
      \def\keyformat##1{\hfil\llap{{##1}\enspace}} 
  \fi\fi\fi\fi}
}
\refstyle{D}   
\def\refsfont@{\tenpoint} 
\def\email{\let\savedef@\email
  \def\email##1\endemail{\let\email\savedef@
  \toks@{\def\usualspace{{\it\enspace}}\endgraf\indent\tenpoint}%
  \toks@@{\tt ##1\par}%
  \expandafter\xdef\csname email\number\addresscount@\endcsname
  {\the\toks@\frills@{{\noexpand\smc E-mail address\noexpand\/}:%
     \noexpand\enspace}\the\toks@@}}%
  \nofrillscheck\email}

\def\D{\text{\rm D}}

\def\GL{\operatorname{GL}}
\def\Aut{\operatorname{Aut}}

\def\Nor{\text{\cal N}}

\def\proof{\noindent{\bf Proof.\ }}


\newcount\secno

\secno=1

\topmatter

\title
Galois theory for a class of modular lattices
\endtitle

\author
Alexandre A.$\,$Panin, Anatoly V.$\,$Yakovlev
\endauthor

\affil
{\it Department of Mathematics and Mechanics\\ 
St.Petersburg State University\\ 
2 Bibliotechnaya square,\\
St.Petersburg 198904, Russia}
\endaffil


\address
\endaddress
\email
alex@ap2707.spb.edu
\endemail
\address
\endaddress
\email
yakovlev@yak.pdmi.ras.ru
\endemail


\abstract
We construct Galois theory for sublattices of certain
complete modular lattices and their automorphism groups.
A well--known description of the intermediate subgroups
of the general linear group over an Artinian ring containing
the group of diagonal matrices, due to Z.I.Borewicz and
N.A.Vavilov, can be obtained as a consequence of this theory.
Bibliography: 11 titles.
\endabstract

\endtopmatter

\document

\heading
\S~\the\secno. Introduction
\endheading
\advance\secno by 1

The description of subgroups in the general linear group
over a semilocal ring $R$ containing the group of diagonal matrices
was obtained in the series of papers of Z.I.Borewicz and N.A.Vavilov
[Bo2], [BV], [V1], [V2]. One may find a wealth of background information
and many further related references in the surveys [V3], [V4].

A.Z.Simonian [S] showed how to state these results in terms of Galois
correspondences between sublattices of certain lattices and subgroups
of their automorphism groups for the case, when $R$ is a field.

Galois theory for lattices is constructed in the present paper. The
description of the intermediate subgroups in the general linear
group over an Artinian ring, containing the group of diagonal matrices,
can be deduced from the results proved here.

\smallskip

Let $L$ be a lattice and $G$ a subgroup of the group $\Aut(L)$ of all
automorphisms of the lattice $L$. Consider a subgroup $F$ of the
group $G$ and a sublattice $M$ of the lattice $L$. By definition, put

\smallskip

$L(F)=\{l\in L$ such that $f(l)=l$ for every $f\in F\}$,

\smallskip

$G(M)=\{g\in G$ such that $g(m)=m$ for every $m\in M\}$

\smallskip\noindent
(it is clear that $L(F)$ is a sublattice of $L$ and $G(M)$ is a subgroup
of $G$).

\smallskip

Let $L_0$ be a sublattice of $L$; we put $H=G(L_0),\ L_0^{'}=L(H)$.
The set of sublattices of $L_0^{'}$ is denoted by {\eufrak M} 
and the set of subgroups of $G$ containing $H$ by {\eufrak N}.
We define mappings $\varphi:$ {\eufrak M} $\to$ {\eufrak N} as
$\varphi(M)=G(M)$ for $M\in$ {\eufrak M} and
$\psi:$ {\eufrak N} $\to$ {\eufrak M} as
$\psi(F)=L_0^{'}(F)$ for $F\in$ {\eufrak N}. It's easy to see that
$(\varphi,\psi)$ is a Galois correspondence between {\eufrak M} and
{\eufrak N}.

\smallskip

We denote hereafter the operations ``infimum'' and ``supremum'' in an arbitrary
lattice as $\cdot$ and $+$, correspondingly.

If $M$ is a lattice, $x_1,\ldots,x_s\in M$, then for every
$i,\ 1\leqslant i\leqslant s$, we put
$\widehat{x}_i=x_1+\ldots+x_{i-1}+x_{i+1}+\ldots+x_s$.

\heading
\S~\the\secno. Formulation of the main result
\endheading
\advance\secno by 1

Let $L$ be a modular lattice of finite length, $L_0$ its sublattice, which is
a Boolean algebra, $G$ a subgroup of the group of all automorphisms
of the lattice $L$, $H=G(L_0)$.

Let $e_1,e_2,\ldots,e_n$ be the atoms of $L_0$, $d(-)$ the dimension
function on the lattice $L$. We require that the following conditions are
fulfilled (it is supposed that, unless otherwise stated, the values of
all indices are changing from $1$ to $n$):

\smallskip

$1^0.\ \ 0_{L_0}=0_L,\ 1_{L_0}=1_L$.

\smallskip

$2^0.$ The function $d$ is constant on the set of atoms of $L_0$;
we denote its value by $m$.

\smallskip

If there are two collections of elements in $L$, namely, $(x_1,\ldots,x_k)$
and $(y_1,\ldots,\break y_k)$, then we write
$(x_1,\ldots,x_k)\leqslant (y_1,\ldots,y_k)$, if $x_i\leqslant y_i$
for every $i,\ 1\leqslant i\leqslant k$. We define also the ``infimum'' and
``supremum'' of two such collections coordinatewise.

For every $x\in L$ its support $[x]$ is defined as the minimal (with respect
to the ordering introduced above) collection $(x_1,\ldots,x_n)$, where

\smallskip

\quad\quad\quad $x_i\leqslant e_i\ \ and\ \ x\leqslant x_1+\ldots+x_n$
\hfill $(+)$

\smallskip

It is proved in \S$\,$4 that the support is well defined.
We put $[x]=([x]_1,\ldots,[x]_n)$.

Let's denote by $H_i$ the set of automorphisms of $H$ which do not change
all elements $x\in L$ such that $[x]_i=0$.

For every $i\neq j$ and every $x\leqslant e_j$
we denote by $H_{ij}(x)$ the set of $f\in G$ such that:

\smallskip

$1)\ f(x_s)=x_s$ for every $s\ne i,\ x_s\leqslant e_s$

\smallskip

$2)\ [f(x_i)]_i=x_i$ for every $x_i\leqslant e_i$

\smallskip

$3)\ [f(e_i)]_k=\cases 0,&k\neq i,j;\\ e_i,&k=i;\\ x,&k=j\\ \endcases$

\smallskip

\noindent (note that $H_{ij}(x)$ may be empty). Elements of
$H_{ij}(x)$ will be called transvections.

We denote by $\overline{L}_0$ the set of elements of the form
$\sum\limits_{i=1}^{n}x_i$, where $x_i\leqslant\ e_i$.
It will be proved (see \S$\,$4) that it follows from the already imposed
conditions that $\overline{L}_0$ is a sublattice of $L$.

We require that the following additional conditions are fulfilled:

\smallskip

$3^0$. If $a\in G$ and $[a(e_i)]_i=e_i$ for some $i$, then there exists
$h\in H_i$ such that $[ha(x_i)]_i=[ah(x_i)]_i=x_i$ for every
$x_i\leqslant e_i$.

\smallskip

$4^0$. There exists $h\in H_t\cap G(\overline{L}_0)$ such that
$[aha^{-1}(x_i)]_r = [a([a^{-1}(x_i)]_t)]_r$ for every
$a\in G,\ r\neq i,\ x_i\leqslant e_i$.

\smallskip

$5^0$. Let $u\in\overline{L}_0,\ u\geqslant e_i$ for some
$i;\ g\in G,\ [g(u)]_i=e_i$. Then there exists $t\in G$ such that:

$1)\ [gt(e_i)]_i=e_i$,

$2)\ t(e_s)=e_s$ for every $s\neq i$,

$3)\ [t(e_i)]_j\leqslant [u]_j$.

\smallskip

$6^0$. If $f,g\in G$ are such that $[f(e_i)]_j\leqslant [g(e_i)]_j$
for some $i,j$, then $[f(x)]_j\leqslant [g(x)]_j$ for every
$x\in L_0^{'},\ x\leqslant e_i$.

\smallskip

$7^0$. If $u\leqslant e_j$ for some $j$, then for every $i\neq j$
there exist $y_1\leqslant e_j,\ldots,y_s\leqslant e_j$ such that
$u=\sum\limits_{r=1}^{s}y_r$ and $H_{ij}(y_r)\neq\varnothing$.

\smallskip

$8^0$. If $x=[f(e_i)]_j$ for some $f\in G,\ i\neq j$, then there exists
$g\in H_{ij}(x)$ such that $[g(u)]_j=[f(u)]_j$ for every $u\leqslant\ e_i$.

\smallskip

$9^0$. If $w\in L,\ d(w)=m,\ $ and $[w]=(0,\ldots,e_i,\ldots,x,\ldots,0)$,
where $H_{ij}(x)\neq\varnothing$, then there exists
$t\in H_{ij}(x)$ such that $t(w)=e_i$.

\smallskip

$10^0$. Let $a_1\in H_{ij}(x_1),\ldots,a_s\in H_{ij}(x_s)$ and
$y\leqslant x_1+\ldots+x_s$ be such that $H_{ij}(y)\neq\varnothing$. Then
$H_{ij}(y)\subseteq \langle H,a_1,\ldots,a_s\rangle$.

\smallskip

$11^0$. If $a\in G$, then for every $t,\ i\neq j$ and every $h\in H_t$ the set
$H_{ij}([aha^{-1}(e_i)]_j)\cap\break \langle a,H\rangle$ is not empty.

\proclaim
{Theorem 2.1}
For every subgroup
$F\geqslant H$ of the group $G$ there exists a sublattice $K$ of $L_0^{'}$
such that $G(K)\trianglelefteq F$. Moreover, if it is assumed that
the lattice $\overline{L}_0(H)$ is finite, then the index $(F:G(K))$ is
finite.
\endproclaim

\noindent {\bf Remark.} The lattice $K$ is not uniquely determined (see
\S$\,$10).

\smallskip

\S \S$\,$4--8 are devoted to the proof of Theorem 2.1. The case $n=1$
is trivial, therefore we assume hereafter that $n\geqslant 2$.

\heading
\S~\the\secno. The case $m=1$
\endheading
\advance\secno by 1

A.Z.Simonian [S] investigated the Galois correspondence introduced in \S$\,$1
for $m=1$.

Theorem 2.1 and results of \S$\,$10 on uniqueness imply

\proclaim
{Theorem 3.1}
Let $L$ be a modular lattice of finite length,
$L_0$ its sublattice of the same length which is a Boolean algebra,
$G$ a subgroup of the group of all automorphisms of the lattice
$L,\ H=G(L_0),\ L_0=L(H)$.
Provided that the conditions $1^{'}-4^{'}$ stated below are
fulfilled, for every subgroup $F\geqslant H$ of the group $G$
there exists a unique sublattice $K$ of the lattice $L_0$ containing $0$
and $1$ such that $G(K)\trianglelefteq F$ and $(F:G(K))<\infty$.
\smallskip\indent
$1^{'}$. There exists at least one automorphism from $H_i$,
which changes all atoms $x\in L \setminus \{ e_i \} $ such that
$[x]_i=e_i$.
\smallskip\indent
$2^{'}$. If $x,y\in L$ are atoms with $[x]=[y]$, then there exists
$h\in H$ such that $h(x)=y$.
\smallskip\indent
$3^{'}$. For every $i\neq j$ the set $H_{ij}(e_j)$ is not empty.
\smallskip\indent
$4^{'}$. If $a\in G$, then for every $\ t,\ i\neq j$ and every
$h\in H_t$ the set $H_{ij}([aha^{-1}(e_i)]_j)\cap \langle a,H\rangle$
is not empty.
\endproclaim

Note that the conditions $1^{'}-4^{'}$ of Theorem 3.1 are not identical
with the conditions of Theorem 2.1 [S], which seem to be simpler than ours.

\heading
\S~\the\secno. Properties of the support
\endheading
\advance\secno by 1

\proclaim
{Lemma 4.1}
If $M$ is an arbitrary modular lattice,
$x,y,z,t\in M$ and $(x+z)\cdot (y+t)=0$, then
$(x+y)\cdot (z+t)=x\cdot z+y\cdot t.$
\endproclaim

\proof See [Bi].

\proclaim
{Corollary 1}
If $x,y\in\overline{L}_0$,
then $x\cdot y=\sum\limits_{i=1}^{n}[x]_i\cdot [y]_i.$
\endproclaim

\proclaim
{Corollary 2}
$\overline{L}_0$ is a lattice.
\endproclaim

If there are two collections satisfying the condition $(+)$ from \S$\,$2, then
so does their ``infimum''. Since there exists at least one collection with the
property $(+)$ (see the condition $1^0$), we see that the
support is well defined.

\proclaim
{Lemma 4.2}
If $M$ is an arbitrary modular lattice,
$x,x_1,\ldots,x_s\in M$, then\break
$\sum\limits_{i=1}^{s}{(x+\widehat{x}_i)\cdot x_i}=
(\sum\limits_{i=1}^{s}x_i)\cdot\prod\limits_{i=1}^{s}{(x+\widehat{x}_i)}.$
\endproclaim

\proof By induction, using the modularity law.

\proclaim
{Lemma 4.3}
For every $x\in L$ and every
$i\ \ [x]_i=(x+\widehat{e}_i)\cdot e_i$.
\endproclaim

\proof By Lemma 4.2
$\sum\limits_{i=1}^{n}{(x+\widehat{e}_i)\cdot e_i}=
\prod\limits_{i=1}^{n}(x+\widehat{e}_i)\geqslant x.$
Further,
$(x+\widehat{e}_i)\cdot e_i\leqslant ([x]_i+\widehat{e}_i)\cdot e_i=[x]_i$,
and we get the desired equality.

\proclaim
{Corollary}
For every $x,y\in L \ [x+y]=[x]+[y]$.
\endproclaim

\proclaim
{Lemma 4.4}
Let $v\in L,\ I\subseteq \{ 1,\ldots,n \}$. Then there
exists $v_I\in L$ such that $[v_I]_i=0$ for
$i\in I,\ v+\sum\limits_{i\in I}[v]_i=v_I+\sum\limits_{i\in I}[v]_i.$
\endproclaim

\proof We put
$v_I=(v+\sum\limits_{i\in I}[v]_i)\cdot (\sum\limits_{i\not\in I}[v]_i)$.
It is clear that
$v_I+\sum\limits_{i\in I}[v]_i\leqslant v+\sum\limits_{i\in I}[v]_i$.
It is easy to verify that the dimensions of the left-hand and right-hand
parts coincide.

\proclaim
{Corollary 1}
For every $t\in H_{ij}(x)$:
\smallskip\indent
$(a)\ \ t(e_i)+e_i=e_i+x$;
\smallskip\indent
$(b)\ \ t(e_i)+x=e_i+x$.
\endproclaim

\proclaim
{Corollary 2}
If $t\in H_{ij}(x)$, then also $t^{-1}\in H_{ij}(x)$.
\endproclaim

\heading
\S~\the\secno. The auxiliary assertions
\endheading
\advance\secno by 1

Hereafter by $F$ we denote a subgroup of the group $G$ containing
the group $H$.

\proclaim
{Lemma 5.1}
For every $a\in F$, indices $t,\ i\neq j$, and every
$h\in H_t$ the set $H_{ij}([aha^{-1}(e_i)]_j)\cap F$ is not empty.
\endproclaim

\proof Follows from the condition $11^0$.

\proclaim
{Lemma 5.2}
Let $u\in L$ be such that
$d(u)=m,\ [u]_i=e_i,\ [u]_j=x$ (for some $i\neq j$),
$H_{ij}(x)\neq\varnothing$. Then there exists $t\in H_{ij}(x)$ such that
$[t(u)]_j=0$.
\endproclaim

\proof Let $i=1,\ j=2$. By Lemma 4.4 there exists
$w\in L$ with the following properties:
$[w]=(e_1,x,0,\ldots,0);\ u\leqslant w+w_1$, where
$w_1=[u]_3+\ldots+[u]_n;\ d(w)=m-d(u\cdot w_1)$.
By Lemma 4.4 $u+\widehat{e}_1=1$. We have $d(u)=m$, therefore
$u\cdot\widehat{e}_1=0$, whence $d(w)=m$. By the condition $9^0$ there exists
$t\in H_{12}(x)$ such that $t(w)=e_1$.
We have $t(u)\leqslant t(w)+t(w_1)\leqslant\widehat{e}_2$, whence
$[t(u)]_2=0$.

\proclaim
{Lemma 5.3}
For every $h\in H$ there exist
$h_i\in H_i,\ i=1,\ldots,n$ such that
$h_nh_{n-1}\ldots h_1h\in G(\overline{L}_0)$.
\endproclaim

\proof One must apply the condition $3^0$.

\smallskip

We define for every $i\neq j$ ``the ideals of transvections''
$$\sigma_{ij}=\sigma_{ij}(F)=
\sum\limits_{x:\ H_{ij}(x)\cap F\neq\varnothing}x\ $$
We also agree that $\sigma_{ii}=e_i$. Note that
$\sigma_{ij}\leqslant e_j$ for every $i,j$.

\proclaim
{Lemma 5.4}
$(i)$ if $H_{ij}(x)\cap F\neq\varnothing$, then
$H_{ij}(x)\subseteq F$.
\smallskip\indent
$(ii)$ for every $y\leqslant\sigma_{ij}$ such that
$H_{ij}(y)\neq\varnothing$ we have $H_{ij}(y)\subseteq F$.
\endproclaim

\proof One must apply the condition $10^0$.

\proclaim
{Lemma 5.5}
$\sigma_{ij}\in L_0^{'}$ for every $i,j$.
\endproclaim

\proof Let $i\neq j,\ h\in H$. For
$L$ is a lattice of finite length, it is sufficient to show that
$h(\sigma_{ij})\leqslant\sigma_{ij}$. Let $f\in H_{ij}(x)\cap F$. Further,
applying the Corollary 1(a) to Lemma 4.4, we obtain $[hf(e_i)]_j=h(x)$.
By the condition $3^0$ it is possible to find $\bar{h}\in H$ such that
$hf\bar{h}\in H_{ij}(h(x))$. Since $hf\bar{h}\in F$, we have
$h(x)\leqslant\sigma_{ij}$, hence $h(\sigma_{ij})\leqslant\sigma_{ij}$.

\smallskip

We denote by $K=K(F)$ the sublattice of the lattice $\overline{L}_0$,
generated by zero and elements $\sum\limits_{j=1}^{n}\sigma_{ij}$, where $i$
changes from $1$ to $n$. By Lemma 5.5 $K$ is a sublattice of $L_0^{'}$.

\heading
\S~\the\secno. Proof of the inclusion $K\subseteq\overline{L_0(F)}$
\endheading
\advance\secno by 1

We denote by $\overline{L_0(F)}$ the lattice which consists of the
elements $l\in\overline{L}_0$ such that $f(l)\in\overline{L}_0$
for every $f\in F$.

\proclaim
{Lemma 6.1}
Let $a\in F$. Let's put
$u_s=\sum\limits_{j=1}^{n}[a(\sigma_{ij})]_s$ for every $s$.
Then $[a^{-1}(u_s)]\leqslant (\sigma_{i1},\ldots,\sigma_{in})$ for every $s$.
\endproclaim

\proof By the definition,
$$[a^{-1}(u_s)]=[\sum\limits_{j=1}^{n}a^{-1}([a(\sigma_{ij})]_s)]=
\sum\limits_{j=1}^{n}[a^{-1}([a(\sigma_{ij})]_s)]$$

Let $j=i$. By the condition $4^0$ there exists $h\in H_s$ such that
$[a^{-1}ha(e_i)]_r=[a^{-1}([a(e_i)]_s)]_r$ for every $r\neq i$.
By Lemma 5.1 we have $[a^{-1}([a(e_i)]_s)]_r\leqslant\sigma_{ir}$
(for every $r$).

Let $j\neq i$. We take an arbitrary $x\leqslant e_j$ such that
$H_{ij}(x)\cap F\neq\varnothing$. Let $b\in H_{ij}(x)\cap F$. Then by the
Corollary 1(a) to Lemma 4.4 $b(e_i)+e_i=e_i+x$, whence

\smallskip

\quad\quad\quad $ab(e_i)+a(e_i)=a(e_i)+a(x)$ \hfill $(\star)$

\smallskip

We show that
$[a^{-1}([ab(e_i)]_s)]+[a^{-1}([a(e_i)]_s)]\leqslant
(\sigma_{i1},\ldots,\sigma_{in})$ for every $s$.

We have already proved that $[a^{-1}([a(e_i)]_s)]_r\leqslant\sigma_{ir}$.
Since $ab\in F$, we have
$[(ab)^{-1}([ab(e_i)]_s)]_r=[b^{-1}a^{-1}([ab(e_i)]_s)]_r\leqslant\sigma_{ir}$.
It is easy to verify that\break $[a^{-1}([ab(e_i)]_s)]_r\leqslant\sigma_{ir}$.

It follows from $(\star)$ that $[a(x)]_s\leqslant [ab(e_i)]_s+[a(e_i)]_s$,
whence $[a^{-1}([a(x)]_s)]\leqslant (\sigma_{i1},\ldots,\sigma_{in})$.

To complete the proof it remains to recall the definition of
$\sigma_{ij}$.

\vfill\eject

\proclaim
{Theorem 6.2}
$\sum\limits_{j=1}^{n}\sigma_{ij}\in\overline{L_0(F)}$ for every $i$.
\endproclaim

\proof Let $a\in F$. We denote
$u=\sum\limits_{j=1}^{n}{\sum\limits_{s=1}^{n}[a(\sigma_{ij})]_s}=
\sum\limits_{s=1}^{n}u_s$, where $u_s=\break
\sum\limits_{j=1}^{n}[a(\sigma_{ij})]_s$.

It is clear that $a(\sum\limits_{j=1}^{n}\sigma_{ij})\leqslant u$.
On the other hand, by Lemma 6.1
$a^{-1}(\sum\limits_{s=1}^{n}u_s)\leqslant \sum\limits_{j=1}^{n}\sigma_{ij}$,
whence $a(\sum\limits_{j=1}^{n}\sigma_{ij})=u\in\overline{L}_0$, hence
$\sum\limits_{j=1}^{n}\sigma_{ij}\in\overline{L_0(F)}$.

\proclaim
{Corollary 1}
$K\subseteq\overline{L_0(F)}$.
\endproclaim

\noindent {\bf Remark.} The lattice $K$ may not coincide with the lattice
$\overline{L_0(F)}$: if $F=H$, then $K=L_0$, but
$\overline{L_0(H)}=\overline{L}_0$.

\proclaim
{Corollary 2}
$G(\overline{L_0(F)})\leqslant G(K)$.
\endproclaim

\heading
\S~\the\secno. Net collections in $L_0^{'}$
\endheading
\advance\secno by 1

\noindent {\bf Definition.} We call by a net collection in $L_0^{'}$
a collection of elements $\tau=(\tau_{ij}),\break i,j=1,\ldots,n$, such that
the following properties are fulfilled (for every $i,j,k$):

\smallskip

$1)\ \ \tau_{ij}\leqslant e_j,$

\smallskip

$2)\ \ \tau_{ii}=e_i,$

\smallskip

$3)\ \ \tau_{ij}\in L_0^{'},$

\smallskip

$4)$ for every $g\in G$ the following assertions are equivalent:

$\quad (i)\ \ [g(e_i)]_j\leqslant \tau_{ij}$

$\quad (ii)\ \ [g(\tau_{ki})]_j\leqslant\tau_{kj}$

\smallskip

\noindent (note that in $4)$ the only nontrivial implications are
$(i)\Rightarrow (ii)$ for the distinct $i,j,k$: see the condition $6^0$).

\proclaim
{Lemma 7.1}
If $\tau^{\alpha}=(\tau_{ij}^{\alpha})$ is a net collection in $L_0^{'}$
for every $\alpha\in I$, then $\tau^{'}=(\tau_{ij}^{'})$, where
$\tau_{ij}^{'}=\prod\limits_{\alpha\in I}\tau_{ij}^{\alpha}$, is also a
net collection in $L_0^{'}$.
\endproclaim

\proof Follows from the definition of a net collection.

\smallskip

For every net collection $\tau=(\tau_{ij})$ in $L_0^{'}$ we denote by
$K_{\tau}$ the sublattice of $L_0^{'}$, generated by zero and elements
$\sum\limits_{j=1}^{n}\tau_{ij},\ i=1,\ldots,n$. Note that the assertion
$4b)$ from the definition of a net collection is equivalent to
$g\in G(K_{\tau})$.

\proclaim
{Theorem 7.2}
$G(K_{\tau})=
\langle H,H_{ij}(x):\ x\leqslant\tau_{ij},\ i\neq j\rangle$.
\endproclaim

\proof We denote
$V=\langle H,H_{ij}(x):\ x\leqslant\tau_{ij},\ i\neq j\rangle$.
By the definition of a net collection $G(K_{\tau})\geqslant V$. Further, let
$g\in G(K_{\tau})$.

A. For every
$i\ \ g(\sum\limits_{j=1}^{n}\tau_{ij})=\sum\limits_{j=1}^{n}\tau_{ij}$,
then $[g(\sum\limits_{j=1}^{n}\tau_{1j})]_1=e_1$. By the condition $5^0$
there exists $t_1\in G$ with the following properties:
$[gt_1(e_1)]_1=e_1;\ t_1(e_k)=e_k$ for
$k\neq 1;\ [t_1(e_1)]_k\leqslant\tau_{1k}$. We show that $t_1\in V$.

By Lemma 5.2 there exists $\bar{t}_2\in H_{12}([t_1(e_1)]_2)\leqslant V$ such
that $[\bar{t}_2t_1(e_1)]_2=0$. Evidently
$[\bar{t}_2t_1(e_1)]_k=[t_1(e_1)]_k$ for $k\neq 2$. Continuing,
we will find $\bar{t}_k\in H_{1k}([t_1(e_1)]_k)\leqslant V,\ k=3,\ldots,n$
such that $\bar{t}_n\ldots\bar{t}_2t_1(e_1)=e_1$, whence
$\bar{t}_n\ldots\bar{t}_2t_1\break\in H$, hence $t_1\in V$.

B. Using the line of reasoning as in the part A, we find $g_1\in V$ such that
$g_1gt_1(e_1)=e_1$. By induction we get $g\in V$.

\proclaim
{Lemma 7.3}
Let $i,j,k$ be pairwise distinct, and
$g\in H_{ij}(x),\ H_{ki}(y)\neq\varnothing$.
Then $H_{kj}([g(y)]_j)\cap \langle H,H_{ij}(x),H_{ki}(y)\rangle
\neq\varnothing$.
\endproclaim

\proof Direct consequence of Lemmas 4.4, 5.2 and 5.3.

\proclaim
{Theorem 7.4}
$\sigma=\sigma(F)$ is a net collection in $L_0^{'}$.
\endproclaim

\proof By Lemma 5.5 $\sigma_{ij}\in L_0^{'}$.
It remains to prove that it follows from $[g(e_i)]_j\leqslant \sigma_{ij}$
that $[g(\sigma_{ki})]_j\leqslant\sigma_{kj}$ (for the distinct $i,j,k$).

Let $y\leqslant\sigma_{ki}$ be such that $H_{ki}(y)\neq\varnothing$.
By the condition $8^0$ we can find $f\in H_{ij}([g(e_i)]_j)\subseteq F$
such that $[f(y)]_j=[g(y)]_j$. Then by Lemma 7.3
$[g(y)]_j\leqslant\sigma_{kj}$.

\smallskip

If $\sigma=\sigma(F)$, then we denote $K=K(F)=K_{\sigma}$.

\proclaim
{Corollary 1}
$G(K)=\langle H,H_{ij}(x)\cap F:\ x\leqslant e_j,\ i\neq j\rangle$.
\endproclaim

\proclaim
{Corollary 2}
$G(K)\leqslant F$.
\endproclaim

\proclaim
{Corollary 3}
$G(\overline{L_0(F)})\leqslant F$.
\endproclaim

\proclaim
{Corollary 4}
The groups $G(K)$ and $F$ have the same transvections.
\endproclaim

\proclaim
{Lemma 7.5}
The groups $G(\overline{L_0(F)})$ and $F$ have the same
transvections.
\endproclaim

\proof Let
$t\in H_{kl}(x)\cap F,\ \sum\limits_{i=1}^{n}x_i\in\overline{L_0(F)}$.
Since $t\in F$, we have
$t(\sum\limits_{i=1}^{n}x_i)=\sum\limits_{i=1}^{n}y_i\in\overline{L}_0$.
It is easy to check that $x_i=y_i$ for $i\neq l$ and $x_l\leqslant y_l$. We
obtain from the equality of dimensions that
$\sum\limits_{i=1}^{n}x_i=\sum\limits_{i=1}^{n}y_i$.

\proclaim
{Lemma 7.6}
If $\tau=(\tau_{ij})$ is a net collection in
$L_0^{'}$, then $\tau=\sigma(G(K_{\tau}))$.
\endproclaim

\proof If $t\in H_{ij}(x)\cap G(K_{\tau})$, then
$x\leqslant \tau_{ij}$. Thus $\sigma_{ij}\leqslant\tau_{ij}$,
but the condition $7^0$ implies that in fact we obtain the equality.

\smallskip

Thus, every net collection in $L_0^{'}$ consists of ``the ideals of
transvections''.

\heading
\S~\the\secno. Proof of the main result
\endheading
\advance\secno by 1

\proclaim
{Lemma 8.1}
$G(\overline{L_0(F)})\trianglelefteq F$.
\endproclaim

\proof Let $g\in G(\overline{L_0(F)}),\ f\in F,\ l\in \overline{L_0(F)}$.
We have $f^{-1}gf(l)=f^{-1}f(l)=l$ (since $f(l)\in\overline{L_0(F)}$),
whence $f^{-1}gf\in G(\overline{L_0(F)})$.

\smallskip

Note that the groups $G(\overline{L_0(F)})$ and $G(K)$ may not coincide
(for $m\neq 1$, see \S$\,$9). If $G(\overline{L_0(F)})\geqslant H$,
then by Lemma 7.5 $G(K)=G(\overline{L_0(F)})$, therefore the relation
$G(K)\trianglelefteq F$ clearly holds true. But general case requires a
special proof.

\proclaim
{Theorem 8.2}
$G(K)\trianglelefteq F.$
\endproclaim

\proof Let $f\in F,\ h\in H_t$. We put
$g=f^{-1}hf$. Since by Lemma 5.1 $[g(e_i)]_j\leqslant\sigma_{ij}$ for every
$i,j$, we obtain $g\in G(K)$ by Theorem 7.4.
Due to Corollary 1 to Theorem~7.4, Lemma 5.3, and the equality
$G(\overline{L}_0)=H_{12}(0)$, for the completion of the proof it
remains to show that for every $f\in F,\ t\in H_{ij}(x)\cap F$,
the automorphism $f^{-1}tf$ belongs to $G(K)$. It follows from Lemmas
7.5 and 8.1 that $t\in G(\overline{L_0(F)})\trianglelefteq F$. Taking into
account the inequality $G(\overline{L_0(F)})\leqslant G(K)$ completes the
proof.

\proclaim
{Lemma 8.3}
If the lattice $\overline{L}_0(H)$ is finite, then the
index of the subgroup $G(K)$ in the group $F$ is finite.
\endproclaim

\proof Since $f^{-1}hf\in G(K)$ for every $f\in F$ and $h\in H$, we have
$hf(\sum\limits_{j=1}^{n}\sigma_{ij})=f(\sum\limits_{j=1}^{n}\sigma_{ij})$.
Hence $f(\sum\limits_{j=1}^{n}\sigma_{ij})\in\overline{L}_0(H)$,
and subject to the settings of Lemma
there is only a finite number of the possibilities for the values
$f(\sum\limits_{j=1}^{n}\sigma_{ij})$. The rest is clear.

\smallskip

Thus, Theorem 2.1 is proved completely.

\proclaim
{Lemma 8.4}
Let $M$ be a sublattice of $\overline{L}_0$.
If $t\in H_{ij}(x)\cap \Nor_GG(M)$, then $t\in G(M)$.
\endproclaim

\proof Let $x\in M$. Then
$x=\sum\limits_{k=1}^{n}x_k,\ x_k\leqslant e_k$.
By the condition $4^0$ there exists $h\in H_j\cap G(\overline{L}_0)$ such
that $[t^{-1}ht(x_i)]_j=[t(x_i)]_j$. By the settings of Lemma
$t^{-1}ht(x)=x$, whence $[t(x_i)]_j\leqslant x_j$. Hence $t(x)=x$.

\heading
\S~\the\secno. Application to linear groups
\endheading
\advance\secno by 1

The description of subgroups in the general linear group $G=\GL(n,R)$
over a semilocal ring $R$, containing the group of diagonal matrices
$D=\D(n,R)$, was obtained in [BV] in terms of nets over $R$.

This description consists in the following. Let $R$ be a semilocal ring
(that is, a ring, quotient of which modulo the Jacobson radical is Artinian),
$C$ its center (which is a commutative semilocal ring). Suppose that all
fields of residues of $C$ modulo its maximal ideals have at least
seven elements. Then for every intermediate subgroup
$F,\ D\leqslant F\leqslant G$, there exists a unique $D$--net
$\sigma$ of order $n$ over $R$ such that
$G(\sigma)\leqslant F\leqslant \Nor(\sigma)$, where $\Nor(\sigma)$ is the
normalizer of $G(\sigma)$ in $G$.

We demonstrate how to deduce this description of the intermediate subgroups 
(frankly, in a slightly weaker form: instead of the semilocal rings we consider
only the Artinian ones) from the results of \S \S$\,$4--8.

\smallskip

Let $R$ be a right Artinian ring, all residue fields of center of which
have at least seven elements, let $V=R^n$ be a free $R$-module of rank
$n$, let $\bar{e}_1=(1,0,\ldots,0),\ldots,\break\bar{e}_n=(0,\ldots,0,1)$ be
the canonical basis of $V$, $e_1=\bar{e}_1R,\ldots,e_n=\bar{e}_nR$.

We denote by $L=L(V)$ the lattice of right submodules of the module $V$
and by $L_0$ the sublattice of the lattice $L$ generated by
$e_1,\ldots,e_n$. It is clear that $L$ is a modular lattice of finite
length, $L_0$ is a Boolean algebra.

Each element $g\in \GL(n,R)$ generates an automorphism of the lattice $L$.
Namely, if $\nu$ is an element of $L$, then
$g(\nu)= \{ g(x):\ x\in\nu \}$. Thus, one can assume that $G=\GL(n,R)$ and
$H=G(L_0)=\D(n,R)$.

Using some elementary facts about the semilocal rings (see [Ba]; [Bo1]; [BV]),
it is easy to verify that the conditions $1^0-10^0$ of Theorem 2.1 hold true.
The assertion of the condition $11^0$ is proved in [BV] (see the proof of
Lemmas 3 and 4).

Thus, we can apply the results obtained in \S \S$\,$4--8.

\smallskip

Let $F$ be an intermediate subgroup, $H\leqslant F\leqslant G$. Then
$G(K)\trianglelefteq F$ for the lattice $K=K_{\sigma}$ (Theorem 8.2).

We put into correspondence to the lattice $K$ the matrix of ideals in
$R\ \ \sigma_K=(\sigma_{ji})$.
Since $\sigma=(\sigma_{ij})$ is a net collection in $L_0^{'}$ (Theorem 7.4),
we see that $\sigma_K$ is a $D$--net of order $n$ over $R$
(see the definition of a net collection; note that the condition $3)$
from this definition implies that $\sigma_{ij}$ is a two-sided ideal
(see the papers cited above)).

Due to the condition $4)$ from the definition of a net collection
the subgroup $G(K)$ coincides with the net subgroup $G(\sigma_K)$.

Thus we proved that for every subgroup $F,\ H\leqslant F\leqslant G$,
there exists a $D$--net $\sigma$ of order $n$ over $R$ such that

\smallskip

\quad\quad\quad $G(\sigma)\leqslant F\leqslant \Nor(\sigma),$\hfill $(\dagger)$

\smallskip

\noindent where by $\Nor(\sigma)$ we denote the normalizer of $G(\sigma)$.

\smallskip

We show that a $D$--net $\sigma$ which satisfies the condition $(\dagger)$
is uniquely determined.

Let $\sigma_1$ and $\sigma_2$ be two $D$--nets, which satisfy
$(\dagger)$. Consider $K_1=K(G(\sigma_1)),K_2=K(G(\sigma_2))$.
Since $G(K_1)=G(\sigma_1),\ G(K_2)=G(\sigma_2)$, we have
$G(K_1)\leqslant F\leqslant \Nor_GG(K_1),
G(K_2)\leqslant F\leqslant \Nor_GG(K_2)$.

Each transvection containing in $G(K_1)$ belongs to $\Nor_GG(K_2)$
and, by Lemma~8.4, is contained in $G(K_2)$. The reverse statement is also
true. Thus, $G(K_1)$ and $G(K_2)$ have the same transvections, whence
$G(\sigma_1)=G(\sigma_2)$, therefore $\sigma_1=\sigma_2$.

Thus we proved that for every subgroup $F,\ H\leqslant F\leqslant G$,
there exists a unique $D$--net $\sigma$ of order $n$ over $R$
such that $G(\sigma)\leqslant F\leqslant \Nor(\sigma)$, hence we obtained
the required description of the intermediate subgroups.

\smallskip

If $R$ is a left Artinian ring, then this description can also be obtained
by passing to the opposite ring.

\smallskip

\noindent {\bf Remark 1.} It follows from Lemma 8.3 that if the lattice
$\overline{L}_0(H)$ is finite, then the index of $G(\sigma)$ in $F$
is finite. The lattice $\overline{L}_0(H)$ consists in our case of
the direct sums of two-sided ideals in $R$.
So, if there is only a finite number of two-sided
ideals in $R$ (for example, if $R$ is a semisimple Artinian ring), then
$(F:G(\sigma))<\infty$. (Indeed, a more powerful result is valid: see [BV]).

\noindent {\bf Remark 2.}  It was mentioned at the beginning of \S$\,$8
that the groups $G(\overline{L_0(F)})$ and $G(K)$ may not coincide.
It is easy to construct examples of such phenomenon for the case of
noncommutative $R$ and $F=H$.

\heading
\S~\the\secno. Appendix
\endheading
\advance\secno by 1

Fix $x\in\overline{L}_0\cap L_0^{'}$ (write $x_k=[x]_k$) and a pair
of indices $i,j$.

We say that $u\leqslant e_j$ satisfies the condition $(\triangle)$, if for
every $f\in G$ it follows from the inequality $[f(e_i)]_j\leqslant u$ that
$[f(x_i)]_j\leqslant x_j$.

\smallskip

\noindent {\bf Example.} $u=0$ satisfies the condition $(\triangle)$ for every
$x,i,j$.

\proclaim
{Lemma 10.1}
If $u_1,u_2$ satisfy the condition $(\triangle)$,
then so does $u_1+u_2$.
\endproclaim

\proof By the condition
$7^0\ \ u_k=\sum\limits_{l=1}^{s_k}y_{kl}$,
where $H_{ij}(y_{kl})\neq\varnothing,\ k=1,2$. By the condition
$8^0$ for every $f\in G$ there exists $g\in H_{ij}([f(e_i)]_j)$ such that
$[g(x_i)]_j=[f(x_i)]_j$, and we can assume that $f\in H_{ij}(y)$ and
$[f(e_i)]_j\leqslant u_1+u_2$.

By the condition
$10^0\ H_{ij}(y)\subseteq \langle H,H_{ij}(y_{kl}),\ k=1,2,\ l=1,\ldots,
s_k\rangle$.

Since $x_s\in L_0^{'}$, we see that $h(x_s)=x_s$ for every $h\in H$ and $s$.

If $t\in H_{ij}(y_{kl})$, then $[t(e_i)]_j=y_{kl}\leqslant u_k$, whence
$t(x_i)\leqslant x_i+x_j$.

\proclaim
{Corollary}
For every $x\in \overline{L}_0\cap L_0^{'}$ and
indices $i,j$ there exists the maximal element
$\tau_{ij}=\tau_{ij}(x)\leqslant e_j$, which satisfies the condition
$(\triangle)$.
\endproclaim

\proclaim
{Lemma 10.2}
$\tau=(\tau_{ij})$ is a net collection in $L_0^{'}$.
\endproclaim

\proof We verify the conditions $1)-4)$ from the definition
of a net collection.

$1)$ Clear.

$2)$ Follows from the condition $6^0$.

$3)$ It is clear that for every $h\in H\ \ h(\tau_{ij})$ satisfies the
condition $(\triangle)$. Then by the definition of
$\tau_{ij}$ we obtain $h(\tau_{ij})\leqslant \tau_{ij}$.

$4)$ By the condition $7^0\ \ \tau_{ki}=\sum\limits_{r=1}^{s}y_r$, where
$H_{ki}(y_r)\neq\varnothing$. Let $g\in H_{ij}(z)$ and $f\in H_{kj}(y)$ be
such that $[g(e_i)]_j\leqslant\tau_{ij}$ and
$[f(e_k)]_j\leqslant\sum\limits_{r=1}^{s}[g(y_r)]_j$.
We have to verify that $[f(x_k)]_j\leqslant x_j$.

By Lemma 7.3 for every $r,\ 1\leqslant r\leqslant s$, there exists
$t_r\in H_{kj}([g(y_r)]_j)\cap \langle H,H_{ij}(z),H_{ki}(y_r)\rangle$.

Since $y_r\leqslant \tau_{ki}$ and $z\leqslant \tau_{ij}$, for every
$t\in H_{ki}(y_r),\ \bar{t}\in H_{ij}(z)$ we have
$t(x_k+x_i+x_j)\leqslant x_k+x_i+x_j$ and
$\bar{t}(x_k+x_i+x_j)\leqslant x_k+x_i+x_j$.
Hence $[t_r(x_k)]_j\leqslant x_j$.

Further, by the condition $10^0\ H_{kj}(y)\subseteq \langle H,t_r,\ r=1,
\ldots,s\rangle$, therefore $[f(x_k)]_j\break\leqslant x_j$. Lemma is proved.

\proclaim
{Lemma 10.3}
If $g\in G$ and $[g(x_i)]_j\leqslant x_j$, then
$[g(e_i)]_j\leqslant\tau_{ij}$.
\endproclaim

\proof If $f\in G$ and $[f(e_i)]_j\leqslant [g(e_i)]_j$,
then by the condition $6^0\ [f(x_i)]_j\leqslant [g(x_i)]_j\leqslant x_j$,
therefore $[g(e_i)]_j$ satisfies the condition $(\triangle)$.

\smallskip

We introduce an equivalence relation on the set of sublattices of
$L_0^{'}$, namely, we put $M_1\sim M_2$, if $G(M_1)=G(M_2)$.

\proclaim
{Theorem 10.4}
Provided that the conditions $1^0-11^0$ of Theorem 2.1
and the following condition
\smallskip\indent
$12^0.\ \ L_0^{'}\subseteq \overline{L}_0$
\smallskip\noindent
are fulfilled, for every subgroup $F,\ H\leqslant F\leqslant G$, there exists
a unique class of the equivalent sublattices of $L_0^{'}$ such that
$G(M)\trianglelefteq F$ for every $M$ of this class.
\endproclaim

\proof Existence follows from Theorem 2.1.

Let $M$ be a sublattice of $L_0^{'}\subseteq\overline{L}_0$. By Lemma 10.2
$\tau(x)=(\tau_{ij}(x))$ is a net collection in $L_0^{'}$ for every
$x\in M$. We put $\tau_{ij}^{'}=\prod\limits_{x\in M}\tau_{ij}(x)$. By
Lemma 7.1 $\tau^{'}=(\tau_{ij}^{'})$ is also a net collection in $L_0^{'}$.
By Theorem 7.2 $G(K_{\tau^{'}})$ is generated by $H$ and its transvections.

We show that $G(K_{\tau^{'}})=G(M)$. Indeed, let $g\in G(K_{\tau^{'}})$.
Then $[g(e_i)]_j\leqslant\tau_{ij}^{'}$ for every $i,j$,
therefore $[g(e_i)]_j\leqslant\tau_{ij}(x)$ for every $x\in M$ and hence
$g\in G(M)$.

If $g\in G(M)$, then $[g(x_i)]_j\leqslant x_j$ for every $x\in M$, and by
Lemma 10.3 $[g(e_i)]_j\leqslant\tau_{ij}(x)$, whence
$[g(e_i)]_j\leqslant\tau_{ij}^{'}$, therefore $g\in G(K_{\tau^{'}})$.

The uniqueness follows now from Lemma 8.4.

\smallskip

It follows from Theorem 10.4 that all closed subgroups in {\eufrak N} (see
\S$\,$1) are of the form $G(K_{\tau})$, where $\tau$ belongs to the
set of net collections in $L_0^{'}$.

Further, the sublattices
$L_0^{'}(\tau)=L_0^{'}(G(K_{\tau}))=\{ x\in L_0^{'}:\ g(x)=x$
for every $g\in G$ such that $[g(e_i)]_j\leqslant\tau_{ij}$ for every $i,j\}$
exhaust all closed sublattices in {\eufrak M}.

Thus, the Galois correspondence introduced in \S$\,$1 is a bijection
between the set of subgroups of the form $G(K_{\tau})$ and the set of
sublattices of the form $L_0^{'}(\tau)$ (both sets are in one-to-one
correspondence with the set of net collections in $L_0^{'}$, see
Lemma 7.6).

\smallskip

We consider now the case $m=1$ in more detail. Note that the condition $12^0$
of Theorem 10.4 is automatically fulfilled in the settings of Theorem 3.1
(since $L_0$ is a Boolean algebra).

It was proved in [S] that all sublattices of $L_0^{'}=L_0$, containing $0$
and $1$, are closed, then each class of the equivalent sublattices of $L_0$,
containing $0$ and $1$, consists of one element, therefore we obtain
uniqueness stated in Theorem 3.1.
Note also that it follows from $G(\overline{L_0(F)})=G(K)$
(see the beginning of \S$\,$8) that $K=\overline{L_0(F)}$.

\smallskip

Let now $L,\ L_0$ and $G$ be as in \S$\,$9. It is easy to prove
that in that case the condition $12^0$ of Theorem 10.4 is fulfilled.
Hence we obtain

\proclaim
{Theorem 10.5}
For every intermediate subgroup
$F,\ \D(n,R)\leqslant F\leqslant \GL(n,R)$, there exists a unique class
of the equivalent sublattices of $L_0^{'}$ such that $G(M)\trianglelefteq F$
for every $M$ of this class.
\endproclaim

It is clear that there is a one-to-one correspondence between net collections
in $L_0^{'}$ and nets of order $n$ over $R$.

Further, closed subgroups in {\eufrak N} are exhausted by the net subgroups,
and closed sublattices in {\eufrak M} by the sublattices
$L_0^{'}(\sigma)=\{ x\in L_0^{'}: \sigma_{ij}x_j\subseteq x_i\},\ \sigma$ is
a net.

\smallskip

It is easy to construct examples showing that a lattice $M$
such that $G(M)\trianglelefteq F$ is not uniquely defined. It is clear that
every class of the equivalent sublattices of $L_0^{'}$ contains
the maximal element, which is a closed sublattice, and the
``canonical sublattice'', generated by zero and sums of ``the ideals of
transvections''. Usually the ``canonical'' one is rather far from its closure,
but it is possible to construct examples showing that this sublattice is
not the minimal element of the class.

\heading
Final remarks
\endheading

$1^o$.~We are sure that Theorem~2.1 can be generalized to a class of complete
modular lattices (not necessary of finite length), and the description
of the intermediate subgroups of the general linear group over a semilocal
ring containing the group of diagonal matrices can be obtained as in \S$\,$9.

\smallskip

$2^o$.~There is a strong analogy between our results and results of N.A.Vavilov
[V5] on the geometry of tori, which is not completely understood at this
moment.

\enddocument